\documentstyle[12pt]{amsart}
\def\cM{{\cal M}}
\def\cMx{{\cal M}_{\xi}}
\def\cU{{\cal U}}
\def\cV{{\cal V}}
\def\cL{{\cal L}}
\def\cO{{\cal O}}

\begin{document}
\baselineskip=15pt

\begin{flushleft}
To appear in {\it Journal f\"ur die reine und angewandte Mathematik}
\end{flushleft}

\title[Poincar\'e bundle and Calabi-Yau hypersurfaces]{Restriction
of the Poincar\'e bundle to a Calabi-Yau hypersurface}

\author{Indranil Biswas}

\address{School of Mathematics, Tata Institute of
Fundamental Research, Homi Bhabha Road, Bombay 400005, India}

\email{indranil@@math.tifr.res.in}

\author{L. Brambila-Paz}

\address{CIMAT, Apdo. Postal 402, C.P. 36240. Guanajuato, M\'exico}

\email{lebp@@fractal.cimat.mx}

\thanks{The second author is a
member of the VBAC Research group of Europroj. The work forms
part of the project Moduli of vector bundles in algebraic
geometry, funded by the EU International Scientific Cooperation
Initiative (Contract no. CI1*-CT93-0031). Both authors
acknowledge the support of DGAPA-UNAM}

\date{}

\maketitle

\section{Introduction}

Let $X$ be a compact connected Riemann surface
of genus $g$, where $g\geq 3$.
Denote by $\cMx \, := \, {\cal M}(n,\xi)$ the moduli space
of stable vector bundles
over $X$ of rank $n$ and fixed determinant $\xi$.
If the degree $\mbox{deg}(\xi)$ and the
rank $n$ are coprime, then there is a universal family
of vector bundles, $\cU$, over
$X$ parametrized by $\cMx$. This family is unique
up to tensoring by a
line bundle that comes from $\cMx$. We fix one
universal family over $X\times  \cMx$
and call it the {\it Poincar\'e bundle}. For any $x \in X$, let
${\cU}_x$ denote the vector bundle over $\cMx$ obtained by
restricting $\cU$ to $x\times \cMx$.
It is known that $\cU$ (see \cite{BBN}) and ${\cU}_x$
(see \cite{NR}
and \cite{Ty}) are stable vector bundles with respect to
any polarization on $X\times \cMx$ and
$\cMx$ respectively.

A smooth anti-canonical
divisor $D$ on $\cMx$ is an example
of a Calabi-Yau variety, i.e., it is connected and
simply connected with trivial canonical line bundle. The
Calabi-Yau varieties are of interest both in string theory and in
algebraic geometry.

In this paper we consider the restrictions of $\cU$ and
${\cU}_x$ to $X\times D$ and $x\times D$ respectively, where
$x\in X$ and $D$ is a smooth anti-canonical
divisor. Denote such restrictions by ${\cU}_D$ and
$({\cU}_D)_x$ respectively.

In Theorem 2.5 and Corollary 2.6 we prove the following~:

\medskip
{\it If $n\geq 3$, then the vector bundle
$({\cU}_D)_x$ is stable
with respect to any polarization on $D$.
Moreover, for the general Riemann surface $X$, the
connected component of the moduli space of semistable
sheaves over $D$, containing the point represented by
$({\cU}_D)_x$, is isomorphic to the Riemann surface $X$.}
\medskip

Actually, we prove that for any point $x$ of any $X$ (not
necessarily the general Riemann surface), the infinitesimal
deformation map for the family ${\cU}_D$ of vector bundles
over $D$ parametrized by $X$, is an
isomorphism from $T_xX$ to $H^1(D,\, \mbox{End}(({\cU}_D)_x))$
[Theorem 2.5]. Therefore, $X$ is an \'etale cover of the
above component of the moduli space of semistable sheaves over $D$.

In Theorems 2.9 and 2.10 we establish the following property
of ${\cU}_D$ when rank $n \geq 3$.

\medskip
{\it For any polarization on
$X\times D$, the vector bundle ${\cU}_D$ is stable.
Moreover, the connected component ${\cal M}^0_X({\cU}_D)$ of the 
moduli space of semistable sheaves
over $X\times D$, containing the point represented
by ${\cU}_D$, is isomorphic to the Jacobian of $X$.}
\medskip

Let $d=\mbox{deg}(\xi)$, and let $X\times \cMx\times 
\mbox{Pic}^0(X) \longrightarrow X\times
{\cal M}(n,d)$ be the map
defined by $(x, E,L) \longmapsto (x, E\bigotimes L)$.

Consider the restriction to $X\times D\times \mbox{Pic}^0(X)$
of the pullback, using the above map, of a universal
vector bundle over $X\times {\cal M}(n, \mbox{deg}(\xi))$.
Theorem 2.10 is proved by establishing that this vector
bundle over $X\times D\times\mbox{Pic}^0(X)$ gives a universal
family of stable vector bundles over $X\times D$ parametrized
by $\mbox{Pic}^0(X)$.

The above results are analogous to those on
${\cU}_x$ and $\cU$ obtained in \cite{NR} and
\cite{BBN} respectively.

It is known that
the connected component, ${\cal M}^0_X(\cU )$, of the moduli
space of semistable sheaves over $X\times \cMx$,
containing the point represented by $\cU$, is isomorphic to
the Jacobian of $X$ \cite{BBN}. As M. S. Narasimhan pointed out
to one of the authors, in order to be able
to recover the Riemann surface $X$ what is still needed is
an isomorphism of ${\cal M}^0_X(\cU )$ 
with $\mbox{Pic}^0(X)$ as polarized varieties.
In the final section we produce a canonical polarization 
on ${\cal M}^0_X(\cU )$ and also on ${\cal M}^0_X({\cU}_D )$.
Using  Torelli's theorem we have (see Theorem 4.3 and Theorem 4.4)

\medskip
{\it Let $X$ and $X' $ be two compact 
connected Riemann surfaces of
genus $g\geq 3$. If ${\cal M}^0_X(\cU ) \cong
{\cal M}^0_{X'}({\cU})$ or
${\cal M}^0_X({\cU}_{D'}) \cong
{\cal M}^0_{X'}(\cU_{D'})$, as polarized varieties,
then $X\cong X'$.}
\medskip

Actually, we give a general construction of a certain line
bundle equipped with a Hermitian structure
on a moduli space, ${\cM}_Y$, of stable
vector bundles over a
compact K\"ahler manifold $Y$ of arbitrary
dimension. The curvature
of the Hermitian connection on this line bundle has
been computed in Theorem 4.1. If $Y = X\times \cMx$ or
$X\times D$, and
${\cM}_Y = \mbox{Pic}^0(X)$, then the curvature form
represents a nonzero multiple of the
natural polarization on $\mbox{Pic}^0(X)$ given by the
cohomology class of a theta divisor.

We note that if the assumption $n \geq 4$ is imposed, then
all the results proved here remain valid for $g =2$.

Restriction of vector bundles to a Calabi-Yau hypersurface
has been considered in \cite{Th}, \cite{DT}.

\medskip
{\it Acknowledgments~:}\, We are grateful to A. N. Tyurin
for pointing out Corollary 2.6. Part of this work was done
during a visit of us to IMATE-MOR, Mich. M\'exico. We wish to
thank Jes\'us Muci\~no-Raymundo and IMATE-MOR for hospitality.

\section{Vector bundles over a smooth anti-canonical 
divisor on the moduli space}

Let $X$ be a compact connected Riemann surface 
of genus $g\geq 3$.
Fix a line bundle $\xi$ over $X$ of degree $d$,
and also fix an integer $n\geq 3$ coprime to $d$.
Let $ {\cMx := \cal M}(n,\xi)$ denote the moduli
space of stable vector
bundles $E$ over $X$ of rank $n$ and with the
determinant ${\bigwedge}^n E$ isomorphic to $ \xi$. 

It is well known that there is a {\it universal vector
bundle} $\cU$ over $X\times \cMx$. A universal vector
bundle $\cU$ is called the {\it Poincar\'e
bundle}. Let $Ad(\cU) \, \subset \, \mbox{End}(\cU )$ be the
subbundle of rank $n^2-1$ defined by the trace zero
endomorphisms. The vector bundle over $\cMx$ obtained
by restricting ${\cU}$ (respectively, $Ad({\cU})$)
to $x\times \cMx$, where
$x\in X$, will be denoted by ${\cU}_x$ 
(respectively, $Ad({\cU})_x$).
 
We know that both the vector bundles
$\cU $ and ${\cal U}_x$
are stable (see \cite{NR} and \cite{BBN}). 
In general, restrictions of stable vector bundles
need not be stable.
In this section we will consider the restrictions 
of ${\cal U}_x$ and $\cU$ to certain subvarieties of
$\cMx$ and $X\times \cMx$ respectively, and prove
that they are stable.

Denoting the canonical line bundle of $\cMx$ by
$K_{\cMx}$, let $D$ be a smooth divisor on $\cMx$ in the 
complete linear system
for the anti-canonical line bundle $K^{-1}_{\cMx}$. Using the
{\it Poincar\'e residue theorem},
the canonical line bundle of $D$ is the 
trivial line bundle. The
variety $D$ is connected since $K^{-1}_{\cMx}$ is ample,
and as $\cMx$ is also simply connected,
$D$ must be simply connected.
In other words, any such divisor $D$ is a
{\it Calabi-Yau variety}.

The restriction of ${\cU}$ to $X\times D$ will be denoted by
${\cU}_D$. For any $x \, \in \, X$, the vector bundle
over $D$ obtained by restricting the
vector bundle ${\cU}_x$ to the subvariety
$x\times D \subset x\times \cMx$,
will be denoted by $({\cU}_D)_x$.
Let $Ad({\cU}_D)$ denote the adjoint bundle of
the vector bundle ${\cU}_D$
over $D$. For any $x \, \in \, X$, denote by
$Ad({\cU}_D)_x$ the vector bundle over $D$
obtained by restricting $Ad({\cU}_x)$ to the subvariety
$x\times D \subset x\times \cMx$.

The divisor $D$ induces the following
exact sequence of sheaves over $\cMx$
$$
0 \, \longrightarrow \, {\cO}_{\cMx}(-D) \,
\longrightarrow \, {\cO}_{\cMx} \, \longrightarrow \,
{\cO}_D \, \longrightarrow \, 0 \, . \eqno(2.1)
$$
Tensoring (2.1) with $Ad({\cU}_x)$ we have the sequence
$$
0 \, \longrightarrow \, Ad({\cU}_x)\otimes {\cO}_{\cMx}(-D)
\, \longrightarrow \, Ad({\cU}_x) \, 
{\buildrel{F}\over\longrightarrow} \,
Ad({\cU}_D)_x \, \longrightarrow \, 0 \, . \eqno(2.2)
$$

The following lemma plays a key r\^ole in the proofs
of the results in this section.

\bigskip
\noindent {\bf Lemma\, 2.1.}\, {\it 
 $H^i(\cMx, \, Ad({\cU}_x)\bigotimes
{\cO}_{\cMx}(-D)) = 0$, for $i = 0,1,2.$}

\bigskip

Lemma 2.1 will be proved in Section 3. Meanwhile,
in this section, Lemma 2.1 will be assumed to be valid.
\bigskip

The restriction homomorphism $F : Ad({\cU}_x)  \longrightarrow
Ad({\cU}_D)_x$ induces a homomorphism
$$
F_i \, : \hspace{.1in} H^i(\cMx,\, Ad({\cU}_x)) \hspace{.1in}
\longrightarrow \hspace{.1in} H^i(D,\, Ad({\cU}_D)_x) \, .
$$ 

We will now deduce from Lemma 2.1 that the two homomorphisms
$F_0$ and $F_1$ are isomorphisms.

\medskip
\noindent {\bf Proposition\, 2.2.}\, {\it The homomorphism
$F_i$ is an isomorphism, for $i =0,1$.}
\bigskip

{\it Proof.}\, This proposition follows by considering the
long exact sequence of cohomology for the exact sequence (2.2)
and then using Lemma 2.1.$\hfill{\Box}$
\medskip

Theorem 2 (page 392) in \cite{NR} says that
$$
H^0(\cMx,\, Ad({\cU}_x)) \, = \, 0
\hspace{.2in} \mbox{and} \hspace{.2in}
H^1(\cMx,\, Ad({\cU}_x)) \, = \, T_x X \, .
$$
Therefore, Proposition 2.2 has the following corollary.

\bigskip
\noindent {\bf Corollary 2.3.}\, {\it For any $x\in X$, 
$$
H^0(D,\, Ad({\cU}_D)_x )\hspace{.1in} = \hspace{.1in} 0\, .
$$
Moreover, $H^1(D,\, Ad({\cU}_D)_x )\cong T_xX$
where $T_x X$ is the holomorphic tangent space of $X$ at $x$.}

\bigskip
Since
$$
H^i(D, \, {\cO}_D) =\left\{\begin{array}{lcl}
{\Bbb C}& \mbox{if} & i=0\\
0& \mbox{if} & i =1 \end{array}\right.
$$
and 
$$
H^i(D, \, \mbox{End}({\cU}_{D})_x) \hspace{.1in} = 
\hspace{.1in} H^i(D,\, {\cO}_D)\oplus H^i(D,\,
Ad({\cU}_D)_x )\, ,
$$
we obtain from Corollary 2.3 that the following two equalities
are valid
$$
H^0(D, \, \mbox{End}({\cU}_{D})_x) \hspace{.1in}
 \cong \hspace{.1in} {\Bbb C}\eqno(2.3)
$$
$$
H^1(D, \, \mbox{End}({\cU}_{D})_x) \hspace{.1in}
\cong \hspace{.1in} T_xX \,.\eqno(2.4)
$$

\bigskip
\noindent {\bf Proposition 2.4.}\,
{\it $\dim H^1(X\times D, \, {\rm End}({\cU}_D)) \,= \,g$.
Moreover, there is a canonical isomorphism $\alpha :\,
H^1(X,{\cal O}) \, \longrightarrow\, H^1(X\times D,
{\rm End}({\cU}_D)).$}
\bigskip

{\it Proof.}\, Let $p:X\times D \longrightarrow X$
be the projection onto the first factor, and let
$$
0  \rightarrow H^1(X,\,{\cal R}_{p*}^0\mbox{End}({\cU}_D))
\rightarrow 
H^1(X\times D, \, \mbox{End}({\cU}_D)) \rightarrow
H^0(X,\, {\cal R}_{p*}^1\mbox{End}({\cU}_D))
\eqno{(2.5)}
$$
be the associated exact sequence obtained from the
Leray spectral sequence. From (2.3) and (2.4) we obtain
the following two canonical isomorphisms
$$
{\cal R}_{p*}^1\mbox{End}({\cU}_D) \hspace{.1in} \cong
\hspace{.1in} TX
$$
and 
$$
{\cal R}_{p*}^0\mbox{End}({\cU}_D) \hspace{.1in} \cong
\hspace{.1in} {\cal O}_X \, .
$$
Now the proposition follows from the exact sequence (2.5).
$\hfill{\Box}$
\bigskip
 
Since $\mbox{Pic}(D) \, \cong  \, {\Bbb Z}$, 
the stability condition of a vector bundle over $D$ does not
depend on the choice of a
polarization. Moreover, the vector space parametrizing
infinitesimal deformations of a vector
bundle $V$ over $D$, namely
$H^1(D, \, \mbox{End}(V))$, coincides with $H^1(D, \,
\mbox{Ad}(V))$.

\bigskip
\noindent {\bf Theorem 2.5.}\, {\it For any
$x\in X$, the vector bundle $({\cU}_D)_x$ over $D$ is stable.
Moreover, the infinitesimal deformation map
$$
T_xX \, \longrightarrow \, H^1(D,\, {\rm Ad}({\cU}_D)_x)
$$
for the family ${\cU}_D$ of vector bundles over $D$
parametrized by $X$, is an isomorphism.}
\bigskip

{\it Proof.}\, We will first prove that $({\cU}_D)_x$
is polystable. The proof of the
polystability of $({\cU}_D)_x$ given below
is similar to the proof of \cite[Proposition 2.4]{BBN}.

For  any $x\, \in \, X$, let
$$
\pi \, : \, Y \hspace{.1in} \longrightarrow \hspace{.1in} X
$$
be a spectral cover which is unramified over $x$. (See \cite{H}
for the definition of a spectral cover; also
can be found in \cite{BNR}.)

The associated rational map
$$
{\overline\pi} \, : \, \mbox{Prym}(Y)  \hspace{.1in}
\longrightarrow \hspace{.1in} \cMx
$$
is a generically finite dominant rational map; $\mbox{Prym}(Y)
\subset \mbox{Pic}^c(Y)$ is the Prym variety for the covering
$\pi$, where $c = \mbox{deg}(\xi) +(g-1)(n^2-1)$.
The subvariety of
$\mbox{Prym}(Y)$ where ${\overline\pi}$ is not defined --
we will
denote this subvariety by $A$ -- is of codimension
at least $3$ and the subvariety of $\cMx$ consisting of
the complement of the image of ${\overline\pi}$ -- we will denote 
this subvariety by $B$ -- is of codimension at least $2$.

The subvariety $D\cap B$ of $D$ is of codimension at least $2$.
Indeed, if it is not true, i.e., $D\cap B$ is a divisor on $D$,
then the fact that $\mbox{Pic}(D) \, \cong \, {\Bbb Z}$ would imply
that $D\cap B$ is actually an ample divisor on 
$D$. This in turn would imply that there is a nonconstant
holomorphic function $f$ on $D \, - \, (D\cap B)$. On the other hand,
since the codimension of $A$ in $\mbox{Prym}(Y)$ is at least $3$, the
pullback of the function $f$ to ${\overline\pi}^{-1}(D- (D\cap B))$,
being
a nonconstant function, would yield the required contradiction.

Using the criterion of \cite[page 7, Lemma 2.1]{BBN}
for a vector bundle to be (semi)stable, expressed in terms of
the extension of a pullback, the following criterion for the
polystability of $({\cU}_D)_x$ is obtained.
In order to prove that the vector
bundle $({\cU}_D)_x$ is polystable, it is enough to show
that $({\overline\pi}\vert_{{\pi}^{-1}(D)})^*({\cU}_D)_x$
on ${\overline\pi}^{-1}(D)$ is polystable with respect to the
restriction of the natural polarization on
$\mbox{Prym}(Y)$ given by the restriction of the
first Chern class of a theta divisor on $\mbox{Pic}^c(Y)$.

Now since $\pi$ is unramified
over $x$, the vector bundle ${\overline\pi}^*({\cU}_D)_x$ 
decomposes as a direct sum of line bundles of same degree. 
More precisely, this vector bundle coincides with the 
restriction of $\bigoplus_{y\in {\pi}^{-1}(x)} L_y$, 
where $L_y$ is the line bundle over $\mbox{Prym}(Y)$ obtained by
restricting a Poincar\'e bundle
over $Y\times \mbox{Pic}^c(Y)$ to $y \times \mbox{Prym}(Y)$.
So $({\overline\pi}\vert_{{\pi}^{-1}(D)})^*({\cU}_D)_x$ must be
polystable, and hence $({\cU}_D)_x$ is polystable.

The Corollary 2.3 says that $({\cU}_D)_x$
is simple. Therefore $({\cU}_D)_x $ is stable.

Let 
$$
T_xX \hspace{.1in} \longrightarrow \hspace{.1in} 
H^1(\cMx \, , {\rm Ad}({\cU})_x) \eqno{(2.6)}
$$
be the infinitesimal deformation map for
the family $\cU$ of vector bundles over $\cMx$ parametrized by
$X$. The infinitesimal deformation map
$$
T_xX \hspace{.1in} \longrightarrow \hspace{.1in}
H^1(D,\, {\rm Ad}({\cU}_D)_x)
$$
for the vector bundle ${\cU}_D$ over $X\times D$,
is simply the composition of the homomorphism in (2.6)
followed by the homomorphism $F_1$ in
Proposition 2.2. In \cite{NR} it has been proved that
infinitesimal deformation map in $(2.6)$ is actually an 
isomorphism. Since from Proposition 2.2
the homomorphism $F_1$ is an isomorphism,
the infinitesimal deformation map
for the family of vector bundles ${\cU}_D$
over $D$, parametrized by
is an isomorphism. This completes the proof
of the theorem.$\hfill{\Box}$
\bigskip

Let ${\cal M}^0(({\cU}_D)_x)$ denote the connected component
of the moduli space of semistable sheaves over $D$ with the
same numerical invariants as $({\cU}_D)_x$ and containing the
point that represents the stable vector bundle $({\cU}_{D})_x$.

\bigskip
\noindent {\bf Corollary \,2.6.}\, {\it For a
general Riemann surface $X$, the moduli space
${\cal M}^0(({\cU}_D)_x)$ is canonically
isomorphic to $X$ .}
\bigskip

{\it Proof.}\, Consider the map
$$
\beta \, : \hspace{.1in}
\hspace{.1in} X \hspace{.1in} \longrightarrow
\hspace{.1in} {\cal M}^0(({\cU}_D)_x)
$$
defined by $x \, \longmapsto \, ({\cU}_D)_x$.
Theorem 2.5 implies that $\beta$ is an \'etale
covering map from $X$ to ${\cal M}^0(({\cU}_D)_x)$. The general 
Riemann surface of genus $g\geq 2$ does not admit a nontrivial 
\'etale covering map to another Riemann surface. By this we
mean that the subset of the moduli space of Riemann surfaces of
genus $g$ representing all Riemann surfaces
having a nontrivial \'etale covering
map to another Riemann surface, is a subvariety of strictly
lower dimension. Thus, for a general Riemann surface $X$,
the map $\beta$ is an isomorphism.$\hfill{\Box}$

\bigskip
\noindent {\bf Remark 2.7.}\, In view of
Corollary 2.6, the moduli space ${\cal M}^0(({\cU}_D)_x)$
constitutes an explicit example of a complete family of
semistable sheaves on a higher dimensional
variety such that all the members of the family
are locally free and also do not have the
numerical invariants of a projectively flat vector bundle.
\bigskip

\noindent {\bf Remark 2.8.}\, Let $X_T \longrightarrow T$ be
a smooth family of irreducible projective curves of genus $g$,
with $g \geq 3$. Suppose we are also given a family of Poincar\'e
bundles, and a family of smooth anti-canonical divisors,
$D_T \longrightarrow T$, parametrized by $T$,
for the family of curves. Let $M_T \longrightarrow T$ denote the
relative moduli space of stable vector bundles for the family
$D_T \longrightarrow T$. Let $p : X_T \longrightarrow M_T$ denote the
morphism which sends any point $x$ in the fiber $X_t$ over $t \in T$
to the stable vector bundle over the anti-canonical divisor
$D_t$ obtained from Theorem 2.5 (by simply restricting the Poincar\'e
bundle to $x \times D_t$). Suppose that $T$ is irreducible, and
the general member in the family of curves $X_T$ does not have
any nontrivial automorphism.
Therefore, over an open subset $U$ of $T$, the map $p$ is an
isomorphism. Since $p$ is a morphism over $T$, it must be an
isomorphism. Consequently, the assertion in
Corollary 2.6 extends to all Riemann surfaces.
\bigskip

We will now consider the vector bundle
$\cU_D$ over $X\times D$.
Actually, the results (and proofs) are
similar to those for ${\cU}$ given in \cite{BBN}.

\bigskip
\noindent {\bf Theorem 2.9.}\, {\it The
restriction of the Poincar\'e bundle
${\cal U}$ to $X\times D$ is stable with
respect to any polarization on $X\times D$.}
\bigskip

{\it Proof.}\, Since $H^1(D,\, {\Bbb Z}) \, =\, 0$,
any polarization $\eta$ on $X\times D$ is the form
$$
\eta \hspace{.1in} = \hspace{.1in} a\lambda \, + \, b\gamma, 
\hspace{.2in} a,b\, >\, 0 \, ,
$$
where $\lambda$ and $\gamma$ are polarizations on
$X$ and $D$ respectively.

Recall Theorem 2.5 which says that the vector bundle
$({\cU}_{D})_x$ is stable for any $x
\in X$. Furthermore, by definition, for any
$\{d \}\in D$, the vector bundle $({\cU}_D)\vert_{X\times \{d\}}$
over $X$ is stable. Hence by Lemma 2.2 in \cite{BBN}, the vector
bundle ${\cU}_D$ is stable with respect to any polarization
on $X\times D$.$\hfill{\Box}$
\bigskip

Let ${\cal M}^0({\cU}_D)$ denote the connected
component of the moduli space of semistable
sheaves over $X\times D$ with the same numerical
invariants as
${\cU}_D$ and containing the point representing
the stable vector bundle ${\cU}_D$.

\bigskip
\noindent {\bf Theorem 2.10.}\,
${\cal M}^0({\cU}_D)\, \cong \, \mbox{Pic}^0(X).$
\bigskip

{\it Proof.}\, Let $p: \, X\times D \, \longrightarrow \,
X$ denote the projection onto the first factor.
Consider the morphism
$$
\delta \, : \hspace{.1in} \mbox{Pic}^0(X) \hspace{.1in}
\longrightarrow \hspace{.1in} {\cal
M}^0({\cU}_D)
$$
defined by $L \, \longmapsto \, {\cU}_D \otimes p^*L.$ 
As in Lemma 3.4 of \cite[page 12]{BBN}, we have that $\delta$ is
injective; we will omit the details. Moreover, 
the Zariski tangent space of ${\cal M}^0({\cU}_D)$ at 
the point $\cU_D \bigotimes p^*L$ is naturally isomorphic
to $H^1(X\times D, \mbox{End}({\cU}_D))$, 
which has dimension $g$ by Proposition 
2.4. From the Zariski's Main Theorem and the fact
that $\mbox{Pic}^0(X)$ is complete we have that $\delta$ is
actually an isomorphism.$\hfill{\Box}$
\medskip

Remark 2.7 is also valid for the moduli
space ${\cal M}^0({\cU}_D)$.

The following section will be devoted to the proof of
Lemma 2.1.

\section{Proof of the main lemma}

As $D$ is linearly equivalent to
$K^{-1}_{\cMx}$, we have 
$$
{\cal O}_{\cMx}(-D) \hspace{.1in} \cong
\hspace{.1in} K_{\cMx}\, . \eqno(3.1)
$$
So Lemma 2.1 is equivalent to the following statement~:
$$
H^i(\cMx , \, Ad({\cU})_x \otimes K_{\cMx}) \, = \, 0
$$
for $i \leq 2$.

The ample generator of $\mbox{Pic}(\cMx) \, \cong \, {\Bbb Z}$
will be denoted by $\Theta$.
We will recall a relationship between the canonical
line bundle and a determinant bundle over $\cMx$.

Define $b \, := \, \mbox{g.c.d}\, (n,d)$, $n' \, := \,
{n\over b}$ and ${\chi}' \, := \, {{d+n(1-g)}\over b}$.

Take a point $y \, \in \, X-x$. Since the 
anti-canonical line bundle
$K^{-1}_{\cMx}$ is isomorphic to $2\Theta$ \cite{Ra},
we conclude that
$$
K_{\cMx} \hspace{.1in} = \hspace{.1in} 
(\mbox{det} (\cU))^{-2n'}\otimes
({\wedge}^n {\cU}_y)^{-2{\chi}'} \, ,\eqno(3.2)
$$
where $\mbox{det} (\cU)$ is the {\it determinant line
bundle} over $\cMx$,
whose fiber over a point represented by a
vector bundle $E$
over $X$ is canonically isomorphic to following line (cf.
\cite{KM})~:
$$
{\bigwedge}^{\mbox{top}}H^0(X,\, E)^*
\bigotimes {\bigwedge}^{\mbox{top}}H^1(X,\, E) \, .
$$
The determinant bundle depends upon the choice of
the universal vector bundle, but the expression of $\Theta$
in terms of the determinant bundle is valid even in
the non-coprime situation where there
is no Poincar\'e bundle.

\medskip
{\it Proof of Lemma 2.1.}\, 
The proof is an application of
the idea initiated in \cite{NR} of computing cohomologies
using the Hecke transformation. We will
not deal with the issues of codimension
computation needed for the application of Hartog type
results on cohomology, since they have
already been resolved in \cite{NR}.

Consider the projective bundle
$$
p \hspace{.1in} :
\hspace{.1in} {\Bbb P}({\cU}_x) \hspace{.1in}
\longrightarrow \hspace{.1in} \cMx
$$
consisting of all hyperplanes in ${\cU}_x$.
Now ${\Bbb P}_x \, := \, {\Bbb P}({\cU}_x)$
parametrizes a natural
family of vector bundles over $X$ of rank 
$n$ and degree $d-1$. Let $\cV$
denote this vector bundle over $X\times 
{\Bbb P}_x$. The restriction of $\cV$ to $X\times \{H\}$,
where $H\subset
E_x$ is a hyperplane representing a point in
${\Bbb P}_x$ over $E
\in \cMx$, fits in the following exact sequence
$$
0 \, \longrightarrow \, {\cV}\vert_{X\times \{H\}} \,
\longrightarrow \, E \, \longrightarrow \, 
{{E_x}\over H} \,
\longrightarrow \, 0 \, ,\eqno{(3.3)}
$$
where $E_x/H$ is supported at $x$.
The general member of the family of vector bundles over $X$,
defined by $\cV$, is stable.

Thus, from a Zariski open subset of $U$ of
${\Bbb P}_x$ there is a natural projection
$$
q \hspace{.1in} : \hspace{.1in} U \hspace{.1in}
\longrightarrow \hspace{.1in} {\cal M}(n, \xi (-x))
$$
to the moduli space of stable vector bundles
of rank $n$ and determinant
$\xi (-x):= \xi\otimes {\cO}_X(-x)$. The key
point in \cite{NR} is that the relative tangent bundle 
over $U$ for the
projection $p$ actually coincides with the 
relative cotangent
bundle for the projection $q$ \cite[page 85]{Ty}.

Considering the exact sequence (3.3), from the definition
of the determinant bundle it follows immediately that
we have the following isomorphism
$$
\mbox{det} (\cV) \hspace{.1in} = \hspace{.1in} p^*\mbox{det} (\cU)
\otimes {\cO}_{{\Bbb P}_x}(1) \, ,\eqno{(3.4)}
$$
where ${\cO}_{{\Bbb P}_x}(1)$ is the tautological line bundle
for the vector bundle ${\cU}_x$.

Let $L$ denote the line bundle over ${\Bbb P}_x$ whose fiber over
the point $(E,H) \in {\Bbb P}_x$ (as in (3.3)) is the kernel of
the homomorphism of the fibers
$$
({\cV}\vert_{X\times \{H\}})_x \hspace{.1in}
\longrightarrow \hspace{.1in} E_x 
$$
in the exact sequence (3.3). It is easy to check that
${\cO}_{{\Bbb P}_x}(1) \, = \, L$. Indeed, this is a consequence of
the exactness of the sequence (3.3) and the fact that the line
bundle ${\bigwedge}^n {\cV}\vert_{x\times {\Bbb P}_x}$
over $x\times{\Bbb P}_x$ is isomorphic to
${\bigwedge}^n {\cU}_x$. The last isomorphism is a consequence
of the fact that the Picard group, $\mbox{Pic}({\Bbb P}_x)$,
of ${\Bbb P}_x$ is discrete and for $z \neq x$, the line bundle
${\bigwedge}^n {\cV}\vert_{z\times {\Bbb P}_x}$ is canonically
isomorphic to ${\bigwedge}^n {\cU}_z$.

Since $p_*T^{\mbox{rel}}_{p} \, = \, Ad({\cU})_x$,
where $T^{\mbox{rel}}_{p}$ is the relative tangent
bundle for the projection $p$, and furthermore the higher 
direct images of
$T^{\mbox{rel}}_{p}$ vanish, from the Leray 
spectral sequence
for the projection $p$ we conclude that the following
$$
H^i(\cMx , \, Ad({\cU})_x \otimes K_{\cMx})\hspace{.1in}
\cong \hspace{.1in}
 H^i(U,\, T^{* ,\mbox{rel}}_{q} \otimes p^*K_{\cMx})
$$
$$
\hspace{.3in} \cong \hspace{.1in} 
H^i(U, \, T^{* ,\mbox{rel}}_{q} \otimes (\mbox{det}
(\cV) \otimes L^* )^{-2n'}\otimes 
({\wedge}^n {\cV}_y)^{-2{\chi}'}) \eqno{(3.5)}
$$
is valid for any $i\geq 0$, where 
$T^{* ,\mbox{rel}}_{q}$ is the relative
cotangent bundle for the projection $q$. Indeed, since $x
\neq y$, the vector bundle ${\cV}_y$ over ${\Bbb P}_x$
is canonically isomorphic to $p^*{\cU}_y$, and
furthermore, as we already
noted earlier, the relative tangent bundle
$T^{\mbox{rel}}_{p}$ is isomorphic to $T^{* ,\mbox{rel}}_{q}$.

Since $K_{\cMx}$ is negative ample, and the restriction of to
a fiber of $q$ is a nonconstant morphism, the 
restriction of $p^*K_{\cMx}$ 
to a fiber of $q$ has strictly negative degree. Furthermore,  
$$
H^i({\Bbb C}{\Bbb P}(r), \, {\Omega}^1(k)) \hspace{.1in}
= \hspace{.1in} 0
$$
if $i\not= r$ and $k <0$ (see \cite{Bo} and
\cite[page 71, Theorem 4.3]{SS}).
Consequently, if $n \geq 4$, then the following equality
$$
{\cal R}^iq_* (T^{* ,\mbox{rel}}_{q} \otimes 
(\mbox{det} (\cV) \otimes L^* )^{-2n'}\otimes 
({\wedge}^n {\cV}_y)^{-2{\chi}'}) \hspace{.1in}=
\hspace{.1in} 0 \eqno{(3.6)}
$$
is valid for $i \leq 2$.

Combining the equality (3.6) with
the Leray spectral sequence applied to the vector
bundle $T^{* ,\mbox{rel}}_{q} \otimes 
(\mbox{det} (\cV) \otimes L^* )^{-2n'}\otimes 
({\wedge}^n {\cV}_y)^{-2{\chi}'}$, for the
projection $q$, we conclude that
$$
H^i(U, \, T^{* ,\mbox{rel}}_{q} \otimes (\mbox{det}
(\cV) \otimes L^* )^{-2n'}\otimes 
({\wedge}^n {\cV}_y)^{-2{\chi}'}) \hspace{.1in} =
\hspace{.1in} 0
$$
for $i\leq 2$ whenevr $n\geq 4$. Now the equality (3.5) yields
that if $n\geq 4$, then
$$
H^i(\cMx , \, Ad({\cU})_x \otimes K_{\cMx})\hspace{.1in} =
\hspace{.1in} 0
$$
for $i\leq 2$. This completes the proof of the lemma when $n \geq 4$.

We now consider the remaining case, namely $n\, = \, 3$.

First note that in this case the left-hand side 
of (3.6) vanishes for $i=0$ and $i=1$. Thus the
proof of lemma will be completed once we are able to
establish the following equality
$$
H^0\big({\cal M}(n, \xi (-x)),\, {\cal R}^2q_* (T^{*
,\mbox{rel}}_{q} \otimes (\mbox{det}
(\cV) \otimes L^* )^{-2n'}\otimes 
({\wedge}^n {\cV}_y)^{-2{\chi}'})\big)
\hspace{.1in} = \, 0 \, . \eqno{(3.7)}
$$

Using a special case of Serre duality, which says that for a
line bundle $\zeta$ over a smooth projective surface $S$
$$
H^2(S,\, {\Omega}^1_S\otimes \zeta) \hspace{.1in} 
= \hspace{.1in}
H^0(S,\, {\Omega}^1_S\otimes {\zeta}^*)^* \, ,
$$
the following isomorphism is obtained
$$
{\cal R}^2q_* (T^{* ,\mbox{rel}}_{q} \otimes 
(\mbox{det}
(\cV) \otimes L^* )^{-2n'}\otimes 
({\wedge}^n {\cV}_y)^{-2{\chi}'}) \ \ \ \ \ \ 
$$
$$
\cong  \hspace{.15in}
\big({\cal R}^0q_* (T^{* ,\mbox{rel}}_{q} \otimes (\mbox{det}
(\cV) \otimes L^* )^{2n'}\otimes 
({\wedge}^n {\cV}_y)^{2{\chi}'})\big)^* \, .
\eqno{(3.8)}
$$

Let $Ad({\overline\cU})$ denote the universal 
adjoint vector bundle over
the smooth locus of $X\times {\cal M}(n, \xi (-x))$. 
We observe that although there may not be any universal
vector bundle over the smooth locus of
$X\times {\cal M}(n, \xi (-x))$, the universal adjoint always
exists. The restriction of $Ad({\overline\cU})$ to $x\times 
{\cal M}(n, \xi (-x))$
will be denoted by $Ad({\overline\cU})_x$. Let
$$
{\overline\Theta} \hspace{.1in} \in \hspace{.1in}
\mbox{Pic}({\cal M}(n, \xi (-x)))
$$
be the positive generator.

Now ${\cal R}^0q_* (T^{* ,\mbox{rel}}_{q} 
\otimes (\mbox{det} (\cV) \otimes L^* )^{2n'}\otimes 
({\wedge}^n {\cV}_y)^{2{\chi}'}) \, \cong \,
Ad({\overline\cU})_x
\bigotimes {\overline\Theta}^{\otimes k}$, where $k >0$. Since
$Ad({\overline\cU})_x$ is semistable of degree
zero, and $Ad({\overline\cU})_x = Ad({\overline\cU})^*_x$, 
we conclude that
$$
H^0\big({\cal M}(n, \xi (-x)), \, \big({\cal R}^0q_* 
(T^{* ,\mbox{rel}}_{q} \otimes
(\mbox{det} (\cV) \otimes L^* )^{2n'}\otimes 
({\wedge}^n {\cV}_y)^{2{\chi}'})\big)^*\big) \hspace{.1in}
= \, 0 \, .
$$

Now (3.8) implies (3.7), and this completes the
proof of the Lemma 2.1.$\hfill{\Box}$
\medskip

If $n\geq 4$, then Lemma 2.1 is also valid for $g=2$.
Consequently, all the results in Section 2 remain valid
for $g=2$ if the condition $n\geq 4$ is
imposed. We note that the special situation where $n=3$ and
$g=2$ has been left out in Theorem 2 of \cite{NR}.

\section{A determinant line bundle over the moduli 
space of vector bundles over a K\"ahler manifold}

Let $Y$ be a compact connected K\"ahler manifold of complex
dimension $d$. Fix cohomology classes $c_i \, \in \, H^{2i}(Y,
\, {\Bbb Q})$, $1\leq i \leq r$. Let $\cM$ be a moduli space of
stable vector bundle $E$ of rank $r$ over $Y$ with $c_i(E) \, = \,
c_i$ and fix the Hilbert polynomial. 
It is known that in general there is no universal vector
bundle over $Y \times \cM$ \cite{Ra}, \cite{Le}.

Fix a point $y_0 \, \in \, Y$. Let $i \, : \, y_0 \,
\longrightarrow \, Y$ be the inclusion of $y_0$ and let $q \, :
\, Y \, \longrightarrow \, y_0$ be the projection. There is a
natural vector bundle $W$ over $Y \times\cM$ of rank $r^2$
with the following
property~: for any $m \in \cM$ if $E$ is the vector bundle over
$Y$ represented by $m$, then the restriction of $W$ to $Y\times
m$ is isomorphic to $\mbox{Hom}((i\circ q)^*E,E)$. The vector
bundle $W$ is constructed using the observation that for any
analytic subset $U \, \subseteq \, \cM$ such
that there is a universal vector bundle over $Y\times U$, any two
universal vector bundles over $Y \times U$ differ by tensoring
with a line bundle pulled back from $U$. This observation is a
simple consequence of the projection formula and the fact that
any automorphism of a stable vector bundle is a scalar
multiplication.

Fix a K\"ahler metric $\omega$ on $Y$. From a celebrated theorem
of Yau and Uhlenbeck \cite{UY} (due to Donaldson for projective
manifolds \cite{Do} and due to Narasimhan and Seshadri \cite{NS}
for Riemann surfaces) for any stable
vector bundle $E$ over $Y$ there is a
unique Hermitian-Einstein connection on $E$. Any two
Hermitian-Einstein metrics differ by multiplication with a
constant scalar. This implies that the above vector bundle $W$
has a natural Hermitian metric induced by
any of the Hermitian-Einstein
metrics over the vector bundles represented in $\cM$.

Given any vector bundle $V$ over $Y \times \cM$, by a general
construction of determinant line bundle \cite{KM}, \cite{BGS}
we have a line bundle $\mbox{det}(V)$ over $\cM$ whose fiber over
$m \in \cM$ is canonically identified with the line
$$
\bigotimes_{i=0}^{d} {\bigwedge}^{\mbox{top}}\big(H^i(Y, \,
V\vert_{Y\times m})\big)^{(-1)^i} \, ,
$$
where $\big({\bigwedge}^{\mbox{top}}F\big)^{-1}$ means
${\bigwedge}^{\mbox{top}}F^*$. This construction of
\cite{KM}, \cite{BGS} gives a homomorphism from the
Grothendieck $K$-group
$K(Y \times \cM)$ to the group of holomorphic line bundles over
$\cM$.

Consider the element
$$
(W \, -\, {\cO}^{\oplus r^2})^{\otimes (d+1)} \hspace{.1in} \in
\hspace{.1in} \, K(Y \times \cM) \, , \eqno{(4.1)}
$$
where ${\cO}^{\oplus r^2}$ is the trivial vector bundle on rank
$r^2$ over $Y \times \cM$. Let
$$
\cL \hspace{.1in} := \hspace{.1in} \mbox{det} \big((W \, -\,
{\cO}^{\oplus r^2})^{\otimes (d+1)}\big) \hspace{.1in}
\longrightarrow \, \cM \eqno{(4.2)}
$$
be the determinant line bundle over $\cM$.

By a general construction of \cite{BGS} (\cite{Q} when $Y$ is a
Riemann surface), using the earlier obtained natural metric on
$W$ and the K\"ahler metric $\omega$ on $Y$ as the input, we have a
Hermitian metric on $\cL$ which is known as the {\it Quillen
metric}. Let $\Omega$ denote the first Chern form for the Hermitian
connection on $\cL$. Our next step will be to identity the form
$\Omega$, and in particular, to calculate the first Chern class of
$\cL$.

Let $\psi \, \in \, H^1(Y,\,{\Bbb Q})\bigotimes
H^1(\cM ,\,{\Bbb Q})$ be
the K\"unneth component of $c_1(W) \, \in \, H^2(Y\times\cM ,\,
{\Bbb Q})$. Denoting the projection of $Y\times\cM$ onto the
$i$-th factor by $p_i$, consider
$$
\gamma \hspace{.1in} := \hspace{.1in} (\psi\cup \psi \cup
p^*_1(c_1)^{d-1})\cap [Y] \hspace{.1in} \in \hspace{.1in}
H^2(\cM , \, {\Bbb Q}) \eqno{(4.3)}
$$
where $c_1$, as defined earlier, is the first Chern class of a
typical vector bundle represented in $\cM$ and $\cap [Y]$ is
the cap product with the oriented generator of $H_{2d}(Y, \,
{\Bbb Q})$.

Let ${\overline\psi}$ be the differential 2-form on $Y\times\cM$
obtained by taking the $(1,1)$-K\"unneth component of the first
Chern form $c_1(W)$ for the Hermitian connection on $W$. Also,
let $\bar c$ denote the (unique) harmonic two form on $Y$
representing the cohomology class $c_1$. Clearly the following
two form
$$
\Gamma \hspace{.1in} := \hspace{.1in} \int_Y{\overline\psi}\cup
{\overline\psi} \cup p^*_1({\bar c})^{d-1}
$$
on $\cM$ is closed, and it represents the cohomology class
$\gamma$ defined in (4.3).

For any $y \, \in \, Y$, let $i_y \, : \, \cM\,
\longrightarrow \, Y\times\cM$ be the inclusion defined by $m
\longmapsto (y,m)$. Denote by $c^{0,2}(W)$ the function from $Y$
to the space of two forms on $\cM$ which sends any $y \in Y$ to
$i^*_yc_1(W)$, the pullback of the first Chern form of $W$.
The volume form ${\omega}^d$ on $Y$, for the K\"ahler form
$\omega$, will denoted by $dV$. Now we are in a position the
describe the curvature form $\Omega$.

\bigskip
\noindent {\bf Theorem \,4.1.}\, {\it The first Chern form
$\Omega$ for the Hermitian connection on $\cL$ has the following
expression:
$$
\Omega \hspace{.1in} = \hspace{.1in} {{d+1}\choose 2}\Gamma \, +
\, (d+1)\int_Y c^{0,2}(W) dV \, .
$$
Furthermore, the first Chern class of $\cL$ is ${{d+1}\choose
2}\gamma$.}
\bigskip

{\it Proof.}\, The main theorem of \cite{BGS} (Theorem
0.1) expresses the first Chern form of $\cL$ in terms of the
Chern character form of $W$ and the Todd forms of the relative
tangent bundle for the projection $p_2$ of $Y\times\cM$ onto
$\cM$. Since the virtual rank of $W  - {\cO}^{\oplus r^2}$ is
zero, the lowest degree term of the Chern character of $(W
- {\cO}^{\oplus r^2})^{\otimes (d+1)}$ is of degree
$2(d+1)$. In other words, if $ch_i(V) \in H^{2i}(Y\times \cM ,
{\Bbb Q})$ denotes the component of the Chern character $ch(V)$
of $V \in K(Y\times \cM)$, then $ch_i\big(W  -
{\cO}^{\oplus r^2}\big) =0$ for $i \leq d$.
Thus there is no contribution of the Todd forms of the
tangent bundle of $Y$ in the
expression of the first Chern form of $\cL$ according to Theorem
0.1 of \cite{BGS}. The first part of the theorem now follows
easily.

The second part of the theorem can be proved using the
Grothendieck-Riemann-Roch theorem. The relevant observation is
that the K\"unneth component of the first Chern class $c_1(W)$
in $H^0(Y,\, {\Bbb Q})\bigotimes H^2(\cM ,\, {\Bbb Q})$
vanishes. The second part of the theorem now follows from the
earlier observation that there is no contribution of the Todd
classes of $Y$ in the Grothendieck-Riemann-Roch formula for the
first Chern class of $\cL$.

The second part of the theorem can also be deduced directly from
the first part. The form $i^*_{y_0}c_1(W)$ on $\cM$ vanishes
identically. Hence by the homotopy invariance of the pullback of
a cohomology class, any $i^*_yc_1(W)$ is an exact form. This
implies that $\int_Y c^{0,2}(W) dV$ is a exact form. Now the
second part of the theorem follows from the first
part.$\hfill{\Box}$
\medskip

If either $H^1(Y, {\Bbb Q}) \, = \, 0$ or $c^{d-1}_1 \, = \, 0$,
then from Theorem 4.1 it follows that $c_1(\cL ) \, = \, 0$. We
will give examples where ${\cL}$ is a nontrivial line bundle.

Let $A$ be an abelian variety of
dimension $g$. For $c \, \in \, H^2(A, \, {\Bbb Z})\cap
H^{1,1}(A)$ let $\mbox{Pic}^c(A)$ denote the component of
$\mbox{Pic}(A)$ consisting of line bundles $L$ with $c_1(L) =c$.
The vector space $H^1(\mbox{Pic}^c(A), \, {\Bbb Q})$ is
canonically identified with $H^1(A, \, {\Bbb Q})^*$.
Indeed, by fixing a point of
$\mbox{Pic}^c(A)$ it gets identified with $\mbox{Pic}^0(A)$
Since
$$
\mbox{Pic}^0(A) \, = \, \mbox{Hom}(H_1(A, \, {\Bbb Z})\, ,
U(1))
$$
the vector space $H^1(\mbox{Pic}^0(A), \, {\Bbb Q})$ is
canonically identified with $H^1(A, \, {\Bbb Q})^*$.
It is known that the K\"unneth component in
$H^1(A,\, {\Bbb Q})\bigotimes
H^1(\mbox{Pic}^c(A), \, {\Bbb Q})$ of the first Chern class of a
universal line bundle over $A\times \mbox{Pic}^c(A)$ is the
above mentioned isomorphism.

Let $a \, \in \, H^2(A, \, {\Bbb Q})^*$ be the Poincar\'e dual of
$c^{d-1} \, \in \, H^{2d-2}(A,\, {\Bbb Q})$.
The earlier mentioned isomorphism between $H^1(\mbox{Pic}^c(A),
\, {\Bbb Q})$ and $H^1(A, \, {\Bbb Q})^*$ identifies
$H^2(A, \, {\Bbb Q})^*$ with
$H^2(\mbox{Pic}^c(A), \, {\Bbb Q})$. Let
$$
{\bar a} \hspace{.1in} \in \hspace{.1in} H^2(\mbox{Pic}^c(A), \,
{\Bbb Q})
$$
be the element corresponding to the element $a$ of
$H^2(A, \, {\Bbb Q})^*$.

Using the earlier remark on the K\"unneth component of the first
Chern class of a universal line bundle over $A\times
\mbox{Pic}^c(A)$, from Theorem 4.1 it is easily deduced that in
this situation the following equality
$$
c_1(\cL) \hspace{.1in} = \hspace{.1in} -\, {{d+1}\choose 2}
{\bar a}
$$
is valid. If $c$ is a principal polarization on $A$ then $\bar
a$ is a principal polarization on $\mbox{Pic}^c(A)$.

If $\omega$ is a translation invariant K\"ahler form on
$A$, then the first Chern
form for the Hermitian connection on $\cL$ coincides with
$-{{d+1}\choose 2}{\overline\omega}$, where $\overline\omega$ is
the unique translation
invariant closed form on $\mbox{Pic}^c(A)$ representing the
cohomology class $\bar a$.

Now we will construct a second example where $\cL$ is
nontrivial. Actually this particular example led us to the
construction of $\cL$.

Let $X$ be a compact connected Riemann surface of genus $g\geq 2$,
Let ${\cMx}$ and ${\cU}$ as before.
Fix a Poincar\'e line bundle ${\cal P} \, \longrightarrow X
\times J$, where $J \, = \, \mbox{Pic}^0(X)$.

It has been proved in \cite{BBN} that
the vector bundle ${\cU}$ is stable
with respect to any polarization on $X\times{\cMx}$.
Furthermore, the vector bundle
$$
p^*_{12}{\cU}\otimes p^*_{13}{\cal P} \hspace{.1in}
\longrightarrow \hspace{.1in} X\times{\cMx}\times J \, ,
$$
where $p_{ij}$ is the projection of $X\times{\cMx}\times J$
onto the product of the $i$-th and the $j$-th factor, is a
complete family of vector bundles over $X\times{\cMx}$. In
other words, $J$ is a component of the moduli space of stable
vector bundles over $X\times{\cMx}$.

Set $Y\, = \, X\times{\cMx}$ and in place of $\cM$ substitute
$J$ realized as a component of the moduli space in the above
fashion. Note that the Jacobian $J$ has a natural polarization
which is defined using the cap product on $H_1(X,\, {\Bbb Q})$.
This polarization on $J$ will be denoted by $\Theta$.

\medskip
\noindent {\bf Proposition \,4.2.}\, {\it The first Chern class
$c_1(\cL)$ coincides with $l\Theta$ where $l$ is an integer not
equal to zero.}
\medskip

{\it Proof.}\, Since $\mbox{Pic}({\cMx}) \, \cong \,
{\Bbb Z}$, any polarization on $X\times {\cMx}$ is of the form
$$
a[X] \, + \, b{\Theta}_{{\cMx}} \, ,
$$
where ${\Theta}_{{\cMx}}$ (respectively, $[X]$) is the positive
generator of the cyclic group $\mbox{Pic}({\cMx})$ (respectively,
$H^2(X,\, {\Bbb Z})$),
and $a$, $b$ are strictly positive integers.

Let $\beta$ denote the K\"unneth component in $H^1(X,\,
{\Bbb Q})\bigotimes H^1(J,\, {\Bbb Q})$ of the first Chern class
$c_1({\cal P}) \, \in \, H^2(X\times J,\, {\Bbb Q})$.

Let $q_2$ denote the projection of $X\times{\cMx}\times J$ onto
the factor ${\cMx}$. For a fixed $x\in X$, the first Chern
class $c_1({\cU}\vert_{x\times {\cMx}}) \, \in \,
H^2({\cMx} , \, {\Bbb Q})$ will be denoted by $\delta$.

We now note that from the second part of Theorem 4.1 the following
equality easily obtained
$$
c_1(\cL) \hspace{.1in} = \hspace{.1in} {{d+1}\choose 2}
\int_{X\times{\cMx}} q^*_2{\delta}^{d-1}\cup
p^*_{12}{\beta}^2 \, , \eqno{(4.4)}
$$
where $d \, = \, \dim X\times {\cMx} \, = \, 1+(r^2-1)(g-1)$
and $\int_{X\times{\cMx}}$ is the Gysin map from
$H^i(X\times{\cMx}\times J, \,
{\Bbb C})$ to $H^{i-2d}(J, \, {\Bbb C})$
defined by integration of forms along $X\times{\cMx}$.

The right-hand side of the equality (4.4) is equal to
$$
{{d+1}\choose 2}
\big(\int_{{\cMx}}{\delta}^{d-1}\big) \int_X{\beta}^2 \, ,
$$
where $\int_X$ is the Gysin map for the projection of $X\times
J$ onto $J$.  Using the isomorphism between $H^2({\cMx} ,\,
{\Bbb Z})$ and $\Bbb Z$, the cohomology class $\delta$ is
$mr+1$, where $m\in {\Bbb Z}$ \cite{Ra}. This implies that the
integer $\int_{{\cMx}}{\delta}^{d-1}$ is nonzero.

Finally, $\int_X{\beta}^2$ coincides with $-\Theta$. This
is because the $(1,1)$-component $\beta$ coincides with
the cohomology class on $X\times J$ defined by the natural
identification
between $H^1(X, \, {\Bbb Q})$ and $H^1(J\, {\Bbb Q})$. This
completes the proof of the proposition.$\hfill{\Box}$
\medskip

Since $\Theta$ is a principal polarization on $J$, it can be
recovered in a unique way from any given integral multiple of
$\Theta$. We equip the moduli space ${\cal M}^0_X(\cU )$
with the polarization defined by the positive one
among $\pm c_1(\cL )$.
Now Proposition 4.2 and the Torelli theorem for
Riemann surfaces combine together to give the
following theorem.

\medskip
\noindent{\bf Theorem 4.3.}\, {\it Let $X$ and $X' $ be two
compact connected Riemann surfaces of genus $g\geq 2$. 
There is an isomorphism between ${\cal M}^0_X(\cU )$ and
${\cal M}^0_{X'}(\cU )$ preserving their polarizations if
and only if $X\cong X'.$}
\medskip

Now we set $ Y=X\times D$, where $D$, as before, is a smooth
anti-canonical divisor on ${\cal M}_{\xi}$, and set ${\cal M}
\, = \, {\cal M}^0({\cal U}_D)$, defined in Section 2.
Assume that $n,g \geq 3$. From
Theorem 2.10 we know that ${\cal M}^0({\cal U}_D) =
\mbox{Pic}^0(X)$.

\bigskip
\noindent {\bf Theorem \,4.4.}\, {\it The first Chern class
$c_1({\cal L})$, where $\cal L$ is the line bundle over
${\cal M}^0({\cal U}_D)$ defined in (4.2), coincides with
a nonzero multiple of $\Theta$ on ${\rm Pic}^0(X)$.
Moreover, let $X$ and $X'$ be two
compact connected Riemann surfaces of
genus $g\geq 3$. If ${\cal M}^0_X(\cU ) \cong 
{\cal M}^0_{X'}(\cU )$, as polarized varieties with
polarizations obtained from $c_1({\cal L})$, then 
$X\cong X'$.}
\bigskip

{\it Proof.}\, We have $H^1(D,\, {\Bbb Q}) \, = \, 0$
and $\mbox{Pic}(D) \, = \, {\Bbb Z}$. The restriction of
the cohomology class $\delta$ (defined in the
proof of Proposition 4.2) on ${\cal M}_{\xi}$ to the
subvariety $D$
is nonzero. Given this situation, the proof of the theorem
is exactly identical to the combination of the proofs of
Proposition 4.2 and Theorem 4.3.$\hfill{\Box}$

%%%%%%%%%%%%%%%%%%%%%%%%%%%%%%%%%%%%%%%%%%%%%%%%%%%%%%%%%%%%%%%%

\end{document}